\documentclass[a4paper, twocolumn]{article}

\usepackage{amsmath}
\usepackage{booktabs}
\usepackage{fullpage}
\usepackage[hidelinks]{hyperref}
\usepackage[UKenglish]{isodate}
\usepackage{pgf-pie}
\usepackage{tcolorbox}
\usepackage{xcolor}
\usepackage{url}
\usepackage{utfsym}

\tcbset{height=0.9cm, width=0.8cm, valign=center, halign=center, left=0cm, right=0cm, colback=white}
\newcommand\card[1]{\begin{tcolorbox}#1\end{tcolorbox}}
\newcommand\emphcard[1]{\begin{tcolorbox}[colback=red!30]#1\end{tcolorbox}}
\newcommand\redpawn{{\LARGE{\textcolor{red}{\usym{2659}}}}}
\newcommand\bluepawn{{\LARGE{\textcolor{blue}{\usym{2659}}}}}

\newcommand\customboard[8]{
  \setlength{\tabcolsep}{0.1cm}
  \begin{tabular}{c c c c}
    #1 & #2 & #3 & #4 \\
    #5 & #6 & #7 & #8 \\
    \customboardmore
}
\newcommand\customboardmore[8]{
    #1 & #2 & #3 & #4 \\
    #5 & #6 & #7 & #8
  \end{tabular}
}
\newcommand\board[8]{
  \setlength{\tabcolsep}{0.1cm}
  \begin{tabular}{c c c c}
    \card{#1} & \card{#2} & \card{#3} & \card{#4} \\
    \card{#5} & \card{#6} & \card{#7} & \card{#8} \\
    \boardmore
}
\newcommand\boardmore[8]{
    \card{#1} & \card{#2} & \card{#3} & \card{#4} \\
    \card{#5} & \card{#6} & \card{#7} & \card{#8}
  \end{tabular}
}

\title{Collapsi is strongly solved}
\author{Michael Young\\University of St Andrews}

\begin{document}

\maketitle

\abstract{
  Collapsi is a two-player game of complete information released in 2025 by
  Mark~S.~Ball of \emph{Riffle Shuffle \& Roll}. Played with two pawns on a
  toroidal board of 16 randomly mixed playing cards, players take it in turns to
  move based on the value of the card they sit on, with the game ending when a
  player has no legal moves.

  The number of possible deals after symmetry breaking is low enough, and the
  game tree shallow enough, to make an exhaustive analysis of the game
  feasible. A solver was written that can find an optimal move for a given board
  position in around 20 milliseconds. A search was applied revealing that the
  first player can force a win in 77.7\% of deals, with the second player able
  to force a win in all others. In 3.5\% of deals the losing player can prolong
  the game to the maximum length of 14 plies; a win can never be forced in fewer
  than 7 plies.
}

\section{The game}

The rules for Collapsi are provided by the designer in an online
document~\cite{rules} and in an online video~\cite{youtube}. In the final
version of the game rules (version 1.3.1) the two players each control one
coloured pawn, with the first player controlling red and the second player
blue. The
board is a set of 16 playing cards taken from a standard deck -- four aces, four
2s, four 3s, two 4s and two jacks -- shuffled and arranged into a $4\times 4$
grid, an example layout being shown in Figure~\ref{fig:board}.The red pawn is
placed on the first revealed jack, the blue pawn on the second revealed jack.
The players then take turns
(\textit{plies}) to move with red going first.

\begin{figure}[ht]
  \centering
  \board A223 4A2{\redpawn} 3A23 {\bluepawn}3A4
  \caption{Example Collapsi board layout}
  \label{fig:board}
\end{figure}

On a player's first turn they are on a jack, and must move their pawn 1 space;
on subsequent turns they must move exactly the number of spaces shown
on the card they begin from, with ace representing 1. Moving is done
orthogonally and the grid is toroidal so that the top edge is joined to the
bottom, and the left to the right. On a given move, a pawn may not enter a given
space twice, and may not end its movement on the space occupied by the
opponent's pawn.

After a player's ply, the card they began from is turned face down and cannot be
entered for the rest of the game. The first player who cannot make a legal move
loses the game, which can therefore last at most 14 plies.

This study considers which player wins in various layouts given perfect play.

\section{Game length}

In analysing games, besides just considering which player has a winning
strategy, it is common to evaluate how many plies a game lasts assuming perfect
play on both sides, as when computing endgame tablebases such as those for
chess~\cite{endgame}.

For Collapsi we consider play where the winning player attempts to win in as few
plies as possible, and the losing player attempts to prolong the game, ideally
to the maximum game length of 14 plies. We refer to this as
\textit{game-length-perfect play}. We can evaluate an end-of-game position to
the number of cards remaining face-up if red wins, and the negative of this
number if blue wins -- this is a function that the red player is trying to
maximise, and the blue player is trying to minimise.

\section{Search space}
\label{sec:search-space}

To find a solution, all possible arrangements of the 16 cards must be
considered, naively a space of $16! = 2.1 \times 10^{13}$ arrangements. However,
various symmetries allow us to reduce this search space by eliminating deals
that are strategically equivalent to others, resulting in a more tractable
number of deals that need to be considered.

First observe that suit plays no part in the game: all 3s are equivalent, all 2s
are equivalent, and so on. The jacks differ only in the colour of the pawns that
are placed on them.

The toroidal nature of the board means that the bottom row can be moved to the
top, or the right column moved to the left edge, without materially affecting
play. Each of these operations can be repeated, so that any row can be chosen as
the top one and any column as the left one. We exploit this by only considering
boards in which the red pawn starts on a jack in the top-left corner, without
loss of generality.

We can also constrain the position of the blue-pawn jack without loss of
generality, to take advantage of reflections in the board. There are 15 possible
positions for the blue-pawn jack, but two such positions are strategically
equivalent if they have the same $(x, y)$ distance between the two jacks,
possibly including an $x$--$y$ swap (diagonal reflection). The
possible distances are shown in Figure~\ref{fig:jack-distances}, with one
representative of each distance highlighted: these 5 positions are the positions
considered for the blue-pawn jack when searching.

\begin{figure}[ht]
  \centering
  \customboard
  {\card{J}} {\emphcard{0,1}} {\emphcard{0,2}} {\card{0,1}}
  {\card{0,1}} {\emphcard{1,1}} {\emphcard{1,2}} {\card{1,1}}
  {\card{0,2}} {\card{1,2}} {\emphcard{2,2}} {\card{1,2}}
  {\card{0,1}} {\card{1,1}} {\card{1,2}} {\card{1,1}}
  \caption{Distances of spaces from top-left jack}
  \label{fig:jack-distances}
\end{figure}

Though these 5 positions cover all possible deals, they are not equally
likely. As Figure~\ref{fig:jack-distances} shows, of the 15 positions for the
second jack: $(0,1)$, $(1,1)$ and $(1,2)$ represent four each; $(0,2)$
represents two; and $(2,2)$ represents only one. Analysis of statistics for all
possible deals should therefore weight these positions accordingly.

Overall, with the position of the red-pawn jack fixed, 5 options for the blue-pawn jack,
and the remaining 14 cards split into three sets of four with the two 4s left
over, the number of deals that must be considered is equal to
\[5~\binom{14}{4} \binom{10}{4} \binom{6}{4} = 15~765~750,\]
a number much more amenable to search than the naive space initially considered.
After weighting the blue-pawn jack positions as described above, we can produce
data for a total of 47~297~250 deals, which correspond proportionally to a
uniformly random deal.

Different deals have game trees of different sizes, but experiments show that
around 100~000 possible games can be played from the start position, on average.

\section{Computation}

A Rust library~\cite{github} was written that can solve any given game position,
and also enumerate all possible deals as described above. Exhaustive experiments
were performed using a 13th Gen Intel Core i5-13500 processor. All 15~765~750
deals can be enumerated in 3.3 seconds and then explored and solved in parallel.

The algorithm used for evaluating game positions was minimax search with
alpha--beta pruning and unlimited depth. This was applied both for the simple
goal of finding a winning move, and for the stricter goal of optimising score
with respect to game length.

Given a particular board, it can be determined whether the current player is in
a winning position, and if so a game-length-perfect move can be computed, in around
0.4 milliseconds. This algorithm is applicable to any game state including those
reached after suboptimal moves have been played, thus satisfying the definition
of a \textit{strongly solved} game~\cite{games-solved}. These very short
computation times suggest that there would be little value in the creation of a
database of game positions, since the results can be computed on demand for most
relevant applications.

The result with game-length-perfect play was determined for all deals in a
search that took 92 minutes using 20 cores. All results can be
reproduced by following the instructions provided in the solver's readme
file~\cite{github}.

\section{Results}

The length of games given game-length-perfect play from both sides is shown in
Table~\ref{tab:game-length}. Since red plays first, red wins games with an odd
number of plies and blue with an even number of plies.

\begin{table}[ht]
  \centering
  \begin{tabular}{r r r c}
    \toprule
    \textbf{Plies} & \textbf{Deals} & \textbf{(\%)} \\
    \midrule
    $\leq 6$ & 0 & 0.0 \\
    7 & 525 & 0.0 \\
    8 & 9~705 & 0.0 \\
    9 & 641~520 & 1.4 \\
    10 & 1~248~372 & 2.6 \\
    11 & 19~830~630 & 41.9 \\
    12 & 7~652~868 & 16.2 \\
    13 & 16~279~812 & 34.4 \\
    14 & 1~633~818 & 3.5 \\
    \bottomrule
  \end{tabular}
  \caption{Length of game with perfect play}
  \label{tab:game-length}
\end{table}

A win can never be forced in as few as 6 plies, and there are only 525 deals (177
up to symmetry) in which red can force a win in the minimum 7 plies, of which
one is shown in Figure~\ref{fig:win-in-7}.

\begin{figure}[ht]
  \centering
  \board {\redpawn}341 1{\bluepawn}12 2123 3234
  \caption{A deal where red wins in 7 plies}
  \label{fig:win-in-7}
\end{figure}

Overall, 36~752~487 deals (77.7\%) are a winning position for red, with the
remaining 10~544~763 (22.3\%) a winning position for blue, a considerable
advantage for the second player, as shown in Figure~\ref{fig:win-chance}.

\begin{figure}[ht]
  \centering
  \begin{tikzpicture}
    \pie[
    radius=2,
    color={red!30, blue!30},
    rotate=90,
    text=pin
    ]{
      77.7/Red,
      22.3/Blue
    }
  \end{tikzpicture}
  \caption{Deals won by each player}
  \label{fig:win-chance}
\end{figure}

\section{Old rules}

Prior to game version 1.2.0, each player was allowed to choose a first move of
1, 2, 3 or 4 spaces from the jack they started on, effectively allowing them to
move to any face-up space on the board. This game had a rather different
strategic status, which was also evaluated by the Rust library. The most notable
difference is the win ratio, with red winning only 37.5\% of deals, compared to
77.7\% for the new rules. The greater range of player choices also results in a
much larger search space, with the total search taking around 5 times as long.

\section{Future work}

Now that perfect moves can be computed easily, a qualitative analysis of perfect
moves might yield strategies that players could follow to improve their play,
especially in the early game where the value of a move can be difficult to
evaluate.

The game rules \cite{rules} also list several variations: \textit{big board},
\textit{shifting board}, \textit{player's choice}, and \textit{solo game}.
These profoundly change the nature of the game, and an analysis of their effect
could be interesting.

\section{Acknowledgements}

Thanks to Prof.~Ian Gent for his research suggestions and comments on an early
version, and to Ben Claydon for some valuable tips on data processing in Rust.

\end{document}